\documentclass[12pt]{article}

\usepackage{amsmath, epsfig, cite}
\usepackage{amssymb}
\usepackage{amsfonts}
\usepackage{latexsym}
\usepackage{graphicx}
\usepackage{enumerate}
\usepackage{indentfirst,graphics,epsfig,psfrag}

\newtheorem{thm}{Theorem}[section]

\newtheorem{lem}[thm]{Lemma}

\newcommand{\qed}{{\hfill\rule{4pt}{7pt}}}
\def\pf{\noindent {\it \textbf{Proof}.} }

\numberwithin{equation}{section}

\makeatletter \@addtoreset{equation}{section} \makeatother

\title {\bf Note on packing of edge-disjoint\\ spanning trees in sparse random
graphs\footnote{Supported by NSFC No. 11071130 and the ``973"
project}}

\author{
{\small Xiaolin Chen, Xueliang Li, Huishu Lian}\\
{\small Center for Combinatorics and LPMC-TJKLC}\\
{\small Nankai University, Tianjin 300071, P.R. China}\\
{\small E-mail: chxlnk@163.com; lxl@nankai.edu.cn; lhs6803@126.com}
   }
\date{}
\begin{document}

\maketitle

\begin{abstract}
The \emph{spanning tree packing number} of a graph $G$ is the
maximum number of edge-disjoint spanning trees contained in $G$. Let
$k\geq 1$ be a fixed integer. Palmer and Spencer proved that in
almost every random graph process, the hitting time for having $k$
edge-disjoint spanning trees equals the hitting time for having
minimum degree $k$. In this paper, we prove that for any $p$ such
that $(\log n+\omega(1))/n\leq p\leq (1.1\log n)/n$, almost surely
the random graph $G(n,p)$ satisfies that the spanning tree packing
number is equal to the minimum degree. Note that this bound for $p$
will allow the minimum degree to be a function of $n$, and in this
sense we improve the result of Palmer and Spencer. Moreover, we also
obtain that for any $p$ such that $p\geq (51\log n)/n$, almost
surely the random graph $G(n,p)$ satisfies that the spanning tree
packing number is less than the minimum degree.\\


\end{abstract}

\section{Introduction}

For a graph $G$ of order $n$, the \emph{spanning tree packing
number}, denoted by $\sigma=\sigma(G)$, is the maximum number of
edge-disjoint spanning trees contained in $G$. The spanning tree
packing problem has long been one of the main motives in graph
theory. In 1961, Nash-Williams \cite{Nash} and Tutte \cite{Tutte}
independently obtained a necessary and sufficient condition for a
graph to have $k$ edge-disjoint spanning trees.
\begin{thm}\cite{Nash, Tutte}\label{NT}
A graph $G=G(V,E)$ contains $k$ edge-disjoint spanning trees if and only if
\begin{equation*}
|E_G(\mathcal{P})|\geq k(|\mathcal{P}|-1)
\end{equation*}
for every partition $\mathcal{P}$ of $V$, where $|\mathcal{P}|$
denotes the the number of sets in $\mathcal{P}$ and
$E_G(\mathcal{P})$ are the crossing edges of $G$, i.e., edges
joining vertices that are in different sets of $\mathcal{P}$.
\end{thm}
In the same papers, they also proved that $\sigma(G)=\lfloor\eta(G)\rfloor$, where
$\eta(G)=\min\limits_{E\subseteq E(G)}{\frac{|E|}{\omega(G-E)-1}}$.

Frieze and Luczak \cite{Frieze} firstly considered the spanning tree
packing number of a random graph and they obtained that for a fixed
integer $k\geq 2$ the random graph $G_{k-\text{out}}$ almost surely
has $k$ edge-disjoint spanning trees. Moreover, Palmer and Spencer
\cite{Palmer1} proved that in almost every random graph process, the
hitting time for having $k$ edge-disjoint spanning trees equals the
hitting time for having minimum degree $k$, for any fixed positive
integer $k$. In other words, considering the random graph $G(n,p)$,
for any fixed positive integer $k$, if $p(n)\leq \frac{\log
n+k\log\log n-\omega(1)}{n}$, the probability that the spanning tree
packing number equals the minimum degree approaches to $1$ as
$n\rightarrow \infty$. Note that for a fixed $k$, $\frac{\log
n+k\log\log n-\omega(1)}{n}$ is the best upper bound for $p(n)$ to
guarantee $\delta(G(n,p))\leq k $ a.s.

On the other hand, in Catlin's paper \cite{Catlin1} it was found
that if the edge probability is rather large, then almost surely the
random graph $G(n,p)$ has $\sigma(G)=\lfloor|E(G)|/(n-1)\rfloor$,
which is less than the minimum degree of $G$. We refer papers
\cite{Catlin1} and \cite{Palmer} to the reader for more details.

A natural question is whether there exists a largest $q(n)$ such
that for every $p\leq q(n)$, almost surely the random graph $G(n,p)$
satisfies that the spanning tree packing number equals the minimum
degree.

In this paper, we partly answer this question by establishing the
following two theorems. The first theorem establishes a lower bound
of $q(n)$ with $q(n)\geq (1.1\log n)/n$. Note that this bound for
$p$ will allow the minimum degree to be a function of $n$, and in
this sense we improve the result of Palmer and Spencer.

\begin{thm}\label{th1}
For any $p$ such that $(\log n+\omega(1))/n\leq p\leq (1.1\log n)/n$,
almost surely the random graph $G\sim G(n,p)$ satisfies that the spanning tree
packing number is equal to the minimum degree, i.e.
\begin{equation*}
\lim_{n\rightarrow \infty}\mathbf{Pr}(\sigma(G)=\delta(G))=1.
\end{equation*}
\end{thm}

The second theorem gives an upper bound of $q(n)$ with $q(n)\leq
(51\log n)/n$.
\begin{thm}\label{th2}
For any $p$ such that $p\geq (51\log n)/n$, almost surely the random graph $G\sim G(n,p)$ satisfies
that the spanning tree packing number is less than the minimum degree, i.e.
\begin{equation*}
\lim_{n\rightarrow \infty}\mathbf{Pr}(\sigma(G)<\delta(G))=1.
\end{equation*}
\end{thm}

The rest of the paper is organized as follows: In Section 2, we list
some basic notations and collect a few auxiliary results. Then we
prove Theorem \ref{th1} in Section 3 and give the proof of Theorem
\ref{th2} in Section 4.

\section{Preliminaries}

\subsection{Notation}

Let $G$ be a graph with vertex set $V(G)$ and edge set $E(G)$. The
number of vertices and edges of $G$ are denoted by $|V(G)|$ and
$|E(G)|$, respectively. Given a set $A\subseteq V(G)$, $\bar{A}$
denotes the set $V(G)\setminus A$, and the subgraph of $G$ induced
by $A$ is denoted by $G[A]$. For two disjoint sets $A,B\subseteq
V(G)$, $E(A,B)$ denotes the set of edges between $A$ and $B$. The
minimum degree of $G$ is denoted by $\delta (G)$ and the maximum
degree by $\Delta (G)$. For more notations we refer to the book
\cite{JU}.

In this paper, we consider the Erd\H{o}s-R\'{e}nyi random graph
$G(n,p)$, which is a graph with $n$ vertices where each of the ${n
\choose 2}$ potential edges appears with probability $p$,
independently. Given a graph property $\mathcal{Q}$, we say that a
random graph $G(n,p)$ has property $\mathcal{Q}$ almost surely
(a.s.), if the probability that the random graph $G(n,p)$ has the
property $\mathcal{Q}$ approaches to 1 as $n\rightarrow\infty$.
Therefore, from now on and throughout the rest of this paper, when
needed we always assume that $n$ is large enough. For a positive
integer $n$ and $0\leq p \leq 1$, we denote by Bin$(n,p)$ the
binomial random variable with parameters $n$ and $p$. $X\sim$
Bin$(n,p)$ means that $X$ and Bin$(n,p)$ have the same distribution.
We always write $\log$ for the natural logarithm.

In this paper, we use the following standard asymptotic notations:
as $n\rightarrow \infty$, $f(n)=o(g(n))$ means that $f(n)/g(n)\rightarrow 0$;
$f(n)=\omega(g(n))$ means that $f(n)/g(n)\rightarrow \infty$;
$f(n)=O(g(n))$ means that there exists a constant $C$ such that $|f(n)|\leq C g(n)$;
$f(n)=\Omega(g(n))$ means that there exists a constant $c>0$ such that $f(n)\geq c g(n)$.

\subsection{Inequalities}

In our proofs, we often use the following inequalities \cite{Noga}.
\begin{lem}(Chernoff's inequality)\label{lem1}
Let $n$ be a positive integer, $p\in[0,1]$ and $X\sim$ $Bin(n,p)$.
For every positive $a$,
\begin{equation*}
\mathbf{Pr}(X<np-a)<\exp\left(-\frac{a^2}{2np}\right)\quad\text{and}\quad \mathbf{Pr}(X>np+a)<\exp\left(-\frac{a^2}{2np}+\frac{a^3}{2(np)^2}\right)
\end{equation*}
In particular, if $a\leq np/2$, then
$$\mathbf{Pr}(X>np+a)<\exp\left(-\frac{a^2}{4np}\right).$$
\end{lem}

\begin{lem}\label{lem2}
For $1\leq k\leq n$,
$$\left(\frac{n}{k}\right)^k\leq {n \choose k}\leq \left(\frac{en}{k}\right)^k .$$
\end{lem}

\section{Proof of Theorem \ref{th1}}

In this section, we first give an upper bound of the minimum degree
in Lemma \ref{lem3}, then we show that for any set $S\subseteq
V(G)$, there are enough edges between $S$ and $\bar{S}$ in Lemmas
\ref{lem5}, \ref{lem6} and \ref{lem7}. Finally, we use these lemmas
and Theorem \ref{NT} to prove  Theorem \ref{th1}.
\begin{lem}\label{lem3}
Let $(\log n+\omega(1))/n\leq p\leq (1.1\log n)/n$ and $G\sim
G(n,p)$. Then $\delta(G)\leq \log n/30$ a.s..
\end{lem}
\pf Let $k=\lfloor\log n/30\rfloor$. It is obvious that if
$\delta(G)\leq k$ a.s., then this is also true for every $p'\leq p$
due to monotonicity. Therefore, it is sufficient to prove that for
$p=(1.1\log n)/(n-k)$, $\delta(G)\leq k$  a.s.

Let $v$ be an arbitrary vertex of $G$. We have
\begin{align*}
\mathbf{Pr}(\text{deg}(v)=k)&=\mathbf{Pr}(\text{Bin}(n-1,p)=k)\\
&= {n-1 \choose k}p^k(1-p)^{n-1-k}\\
&\geq \left(\frac{n-k}{k}\right)^kp^k(1-p)^{n-k}\\
&= (1-o(1))\left(\frac{(n-k)p}{k}\right)^ke^{-p(n-k)}\\
&\geq (1-o(1))(33)^{\log n/30}n^{-1.1}\\
&=\omega(1/n).
\end{align*}
Then we use a basic result in the theory of random graphs due to
Bollob\'{a}s (see e.g. \cite{B2}, Chapter 3) which asserts that if
$\mathbf{Pr}(\text{Bin}(n-1,p)=k)=\omega(1/n)$, then $\delta(G)\leq
k$ a.s. This completes the proof.\qed

A vertex is called {\it small} if its degree is less than or equal
to $\log n/6$, and otherwise it is called {\it large}. Denote by
SMALL and LARGE the set of all small vertices and all large
vertices, respectively. Then we can obtain an important structural
property of random graphs as follows.
\begin{lem}\label{lem4}
If $(\log n+\omega(1))/n\leq p\leq (1.1\log n)/n$, then a.s. the random graph $G\sim G(n,p)$
satisfies the following properties:
\begin{enumerate}[\indent $(1)$]
\item $|\,\text{SMALL}\,|\leq n^{1/2}$;
\item No pair of small vertices are adjacent or share a common neighbor.
\end{enumerate}
\end{lem}
\pf $(1)$ Let $s=\lceil n^{1/2}\rceil$. Assume that there exists a vertex set $S$
with order $s$ such that each vertex $v\in S$ is small, which happens with probability at most
\begin{align*}
&{n \choose s}\left(\sum_{k=0}^{\log n/6}{n-1 \choose k}p^k(1-p)^{n-1-k}\right)^s\\
\leq& \left(\frac{ne}{s}\right)^s\left(\frac{\log n}{6}\left(\frac{6(n-1)e}{\log n}\right)^{\log n/6}p^{\log n/6}e^{-p(n-1-\log n/6)}\right)^s\\
\leq& \left(\frac{ne}{s}\cdot\frac{\log n}{6}\cdot\left(6.6e\right)^{\log n/6}\cdot e^{-\log n+p+(\log n/6)p}\right)^s\\
\leq& \left(\frac{ne}{s}\cdot\frac{\log n}{6}\cdot n^{0.482}\cdot n^{-1}\cdot O(1)\right)^s\\
=&O(n^{-0.01s}).
\end{align*}
It means that a.s. $|\,\text{SMALL}\,|\leq n^{1/2}$.

$(2)$ The probability that $G$ violates property $(2)$ can be bounded as follows:
\begin{align*}
\mathbf{Pr}(G\,\,\text{violates}\,\,(2))&\leq{n \choose 2}\cdot p\cdot\left(\sum_{k=0}^{\log n/6}{n-1 \choose k}p^k(1-p)^{n-1-k}\right)^2\\
\quad&+{n \choose 2}\cdot{n \choose 1}\cdot p^2\cdot\left(\sum_{k=0}^{\log n/6-1}{n-1 \choose k}p^k(1-p)^{n-1-k}\right)^2\\
&=O(n^{-0.01}),
\end{align*}
which implies that property $(2)$ holds a.s. \qed

\begin{lem}\label{lem5}
Let $(\log n+\omega(1))/n\leq p\leq (1.1\log n)/n$ and $G\sim
G(n,p)$. Then a.s. for any vertex subset $S$ such that
$\emptyset\neq S\subseteq$ LARGE and $|S|\leq n/(\log n)^3$,
$|E(S,\bar{S})|\geq (\log n/10)\cdot |S|$ .
\end{lem}
\pf We prove this lemma by contradiction. Assume that there exists a
vertex subset $S$ such that $\emptyset\neq S\subseteq$ LARGE,
$|S|\leq n/(\log n)^3$ and $|E(S,\bar{S})|< (\log n/10)\cdot |S|$.
Then the induced subgraph $G[S]$ contains $|S|$ vertices and at
least $\left(\frac{\log n}{6}|S|-\frac{\log
n}{10}|S|\,\right)/2=\frac{\log n}{30}|S|$ edges. The probability
for the existence of such $S$ can be bounded as follows:
\begin{align*}
&\mathbf{Pr}\Bigg(\bigcup_{|S|\leq n/(\log n)^3 \atop S\subseteq \text{LARGE}}\left(|E(S,\bar{S})|< (\log n/10)\cdot |S|\right)\Bigg)\\
\leq&\sum_{r=1}^{\frac{n}{(\log n)^3}}\left({n \choose r}\cdot \sum_{k=(\log n/30)\cdot r}^{r \choose 2}\left({{r \choose 2} \choose k}p^k\left(1-p\right)^{{r\choose 2}-k}\right)\right)\\
\leq&\sum_{r=1}^{\frac{n}{(\log n)^3}}\left(\left(\frac{ne}{r}\right)^r\cdot{r \choose 2}{\frac{r^2}{2}\choose \frac{\log n}{30}r} p^{\frac{\log n}{30}r}\left(1-p\right)^{{r\choose 2}-\frac{\log n}{30}r}\right)\\
\leq&\sum_{r=1}^{\frac{n}{(\log n)^3}}\left(\left(\frac{ne}{r}\right)^r\cdot\frac{r^2}{2}\left(\frac{15erp}{\log n}\right)^{\frac{\log n}{30}r}e^{-\frac{\log n}{2n}r^2+\frac{(\log n)^2}{30n}r}\right)\\
=&O(n^{-20}),
\end{align*}
which implies the correctness of the lemma. \qed

\begin{lem}\label{lem6}
Let $(\log n+\omega(1))/n\leq p\leq (1.1\log n)/n$ and $G\sim
G(n,p)$. Then a.s. for any vertex subset $S$ such that $n/(\log
n)^3\leq |S|\leq n/2 $, $|E(S,\bar{S})|\geq (\log n/10)\cdot |S|$.
\end{lem}
\pf The Event that there exists a vertex subset $S$ such that
$n/(\log n)^3\leq |S|\leq n/2 $ and $|E(S,\bar{S})|< (\log
n/10)\cdot |S|$ happens with probability at most
\begin{align*}
&\sum_{|S|=n/(\log n)^3}^{n/2}\left({n\choose s}\mathbf{Pr}(|E(S,\bar{S})|< \frac{\log n}{10}|S|)\right)\\
\leq&\sum_{|S|=n/(\log n)^3}^{n/2}\left(\left(\frac{ne}{s}\right)^s\cdot e^{-\frac{1}{2}\left(1-\frac{\log n}{10p(n-s)}\right)^2\cdot(n-s)sp}\right)\\
\leq&\sum_{|S|=n/(\log n)^3}^{n/2}\left(\left(\frac{ne}{s}\right)^s\cdot e^{-\frac{s\log n}{10}}\right)\\
\leq&\sum_{|S|=n/(\log n)^3}^{n/2}\left(\frac{n^{9/10}e}{s}\right)^s\\
\leq& \quad\frac{n}{2}\left(\frac{e(\log n)^3}{n^{1/10}}\right)^{\frac{n}{(\log n)^3}}\\
=&o(n^{-20}),
\end{align*}
which gives precisely what we want.\qed

\begin{lem}\label{lem7}
Let $(\log n+\omega(1))/n\leq p\leq (1.1\log n)/n$ and $G\sim
G(n,p)$. Then a.s. for any vertex subset $S$ such that $1\leq
|S|\leq n/(\log n)^3$, $|E(S,\bar{S})|\geq \delta(G)\cdot |S|$.
\end{lem}
\pf For any set $S\subseteq V(G)$ and $1\leq |S|\leq n/(\log n)^3$.
Let $S=S_1\cup S_2$, where $S_1\subseteq$ LARGE, $S_2\subseteq$
SMALL. Then
$|E(S,\bar{S})|=|E(S_1,\bar{S_1})|+|E(S_2,\bar{S_2})|-2|E(S_1,S_2)|$.
By Lemmas \ref{lem3}, \ref{lem4} and \ref{lem5}, we get that
$|E(S_2,\bar{S_2})|\geq \delta(G) \cdot|S_2|$, $|E(S_1,S_2)|\leq
|S_1|$ and $|E(S_1,\bar{S_1})|\geq (\log n/10) \cdot|S_1|$,
respectively. It follows that
\begin{equation*}
|E(S,\bar{S})|\geq \left(\frac{\log n}{10}-2\right)\cdot |S_1| + \delta(G)\cdot |S_2|\geq \delta(G)\cdot |S|.
\end{equation*}
The proof is thus completed.\qed

At present, we are ready to prove Theorem \ref{th1}.

\noindent\textbf{\emph{Proof of Theorem \ref{th1}:}} Recall that
$(\log n+\omega(1))/n\leq p\leq (1.1\log n)/n$. Consider the random
graph $G\sim G(n,p)$. Obviously, $\sigma(G)\leq\delta(G)$ always
holds. We only need to prove that a.s. $\sigma(G)\geq\delta(G)$. By
Theorem \ref{NT}, it is sufficient to show that for any partition
$\mathcal{P}$ of $V(G)$, $|E_\mathcal{P}(G)|\geq \delta(G)\cdot
(|\mathcal{P}|-1)$.

Given a partition $\mathcal{P}=\{V_1,V_2,\ldots,V_t\}$ with $t\geq
2$. Suppose $|V_1|\geq|V_2|\geq\ldots\geq|V_t|$. We distinguish two
cases to prove the theorem, according to the order of $V_1$.

\noindent\textbf{Case 1.} $|V_1|\geq \frac{n}{2}$.

Since $|V_1|\geq \frac{n}{2}$, then $|\bar{V_1}|\leq \frac{n}{2}$
and $|V_i| \leq \frac{n}{2}$ for $2\leq i\leq t$. By Lemmas
\ref{lem3}, \ref{lem6} and \ref{lem7}, $|E(V_1,\bar{V_1})|\geq
\delta(G)\cdot |\bar{V_1}|$ and $|E(V_i,\bar{V_i})|\geq
\delta(G)\cdot |V_i|$ for $2\leq i\leq t$. Therefore,
\begin{equation*}
|E_\mathcal{P}(G)|=\frac{1}{2}\,\sum_{i=1}^{t}|E(V_i,\bar{V_i})|\geq \frac{1}{2}\,\delta(G)\cdot|\bar{V_1}| +\frac{1}{2}\,\delta(G)\cdot\sum_{i=2}^{t}|V_i|.
\end{equation*}
Note that $|\bar{V_1}|=\sum_{i=2}^t|V_i|\geq t-1$. We can conclude that
\begin{equation*}
|E_\mathcal{P}(G)|\geq \frac{1}{2}\,\delta(G)\cdot (t-1)+\frac{1}{2}\,\delta(G)\cdot (t-1)=\delta(G)\cdot(t-1).
\end{equation*}
\textbf{Case 2.} $|V_1|< \frac{n}{2}$.

In this case, we consider two subcases, according to the value of
$t$.

\textbf{Subcase 2.1.} $t\geq2n^{\frac{1}{2}}$.

Let $\mathcal{P}_1=\{V_i|\,1\leq i\leq t, V_i \text{ contains no
small vertex}\}$ and
$\mathcal{P}_2=\mathcal{P}\setminus\mathcal{P}_1$. Then we have that
\begin{equation*}
|E_\mathcal{P}(G)|=\frac{1}{2}\sum_{i=1}^{t}|E(V_i,\bar{V_i})|=\frac{1}{2}\left(\sum_{V_i\in \mathcal{P}_1}|E(V_i,\bar{V_i})|+\sum_{V_j\in \mathcal{P}_2}|E(V_j,\bar{V_j})|\right).
\end{equation*}
Note that $\delta(G)\leq \log n/30$ a.s. For any $V_i\in
\mathcal{P}_1$, by Lemmas \ref{lem5} and \ref{lem6},
$|E(V_i,\bar{V_i})|\geq (\log n/10)\cdot|V_i|$ and for any $V_j\in
\mathcal{P}_2$, by Lemmas \ref{lem6} and \ref{lem7},
$|E(V_j,\bar{V_j})|\geq \delta(G)\cdot|V_j|$. Moreover, by Lemma
\ref{lem4}, $|\mathcal{P}_2|\leq n^{1/2}$. Therefore,
\begin{align*}
|E_\mathcal{P}(G)|&\geq \frac{1}{2}\left(\frac{\log n}{10}\cdot |\mathcal{P}_1|+\delta(G)\cdot |\mathcal{P}_2|\right)\\
&\geq \frac{1}{2}\left(3\,\delta(G)\cdot \left(t-n^{\frac{1}{2}}\right)+\delta(G)\cdot n^{\frac{1}{2}}\right)\\
&\geq \frac{1}{2}\left(3\,\delta(G)\cdot t-2\delta(G)\cdot n^{\frac{1}{2}}\right)\\
&\geq \delta(G)\cdot t>\delta(G)\cdot(t-1).
\end{align*}

\textbf{Subcase 2.2.} $t<2n^{\frac{1}{2}}$.

Note that by Lemmas \ref{lem3}, \ref{lem6} and \ref{lem7}, for any
$1\leq i\leq t$, $|E(V_i,\bar{V_i})|\geq \delta(G)\cdot |V_i|$ a.s.
Then we get that
\begin{align*}
|E_\mathcal{P}(G)|&=\,\frac{1}{2}\sum_{i=1}^{t}|E(V_i,\bar{V_i})|\,\geq\, \frac{1}{2}\delta(G)\cdot\sum_{i=1}^{t}|V_i|\\
&=\,\frac{1}{2}\delta(G)\cdot n \,>\,2\delta(G)\cdot n^{\frac{1}{2}}\,>\,\delta(G)\cdot(t-1).
\end{align*}
Combining the two cases discussed above, we can conclude that $G$
has $\delta(G)$ edge-disjoint spanning trees. It immediately implies
that $\sigma(G)\geq\delta(G)$. We thus complete the proof of Theorem
\ref{th1}.\qed

\section{Proof of Theorem \ref{th2}}

\noindent {\bf\emph{Proof of Theorem \ref{th2}:}} Recall that $G\sim
G_(n,p)$ with $p\geq 51\log n/n$. We first bound the minimum degree
and the maximum degree of $G$. Let $v$ be an arbitrary vertex of
$G$. Then deg($v$), the degree of $v$, obeys the binomial
distribution $\text{Bin}(n-1,p)$. By
$\mathbf{E}(\text{deg}(v))=(n-1)p$ and Chernoff's inequality
\ref{lem1},
\begin{equation*}
\mathbf{Pr}\left(\text{deg}(v)\geq\frac{3}{2}(n-1)p\right)\leq \exp\left(-\frac{(n-1)p}{16}\right)=o(n^{-2}).
\end{equation*}
Hence, by the union bound, with probability at least $1-o(n^{-1})$, $\Delta(G)\leq \frac{3}{2}(n-1)p$.

On the other hand,
\begin{equation*}
\mathbf{Pr}\left(\text{deg}(v)\leq\frac{4}{5}(n-1)p\right)\leq \exp\left(-\frac{(n-1)p}{50}\right)=o(n^{-1.01}).
\end{equation*}
By the union bound again, it follows that with probability at least $1-n^{-0.01}$, $\delta(G)\geq \frac{4}{5}(n-1)p$.
Then we can deduce that a.s.
\begin{equation*}
\sigma(G)\leq \frac{|E(G)|}{n-1}\leq \frac{\Delta(G)\cdot n}{2(n-1)}\leq \frac{3}{4}np< \frac{4}{5}(n-1)p\leq \delta(G),
\end{equation*}
the proof is thus completed.\qed

\end{document}